\relax
\def\DATE{\relax}
\magnification=1100
\baselineskip=14truept
\voffset=.75in
\hoffset=1truein
\hsize=4.5truein
\newdimen\hsizeGlobal
\hsizeGlobal=4.5truein%
\vsize=7.75truein
\parindent=.166666in
\pretolerance=500 \tolerance=1000 \brokenpenalty=5000

\footline={\vbox{\hsize=\hsizeGlobal\hfill{\rm\the\pageno}\hfill\llap{\sevenrm\DATE}}\hss}

\def\note#1{%
  \hfuzz=50pt%
  \vadjust{%
    \setbox1=\vtop{%
      \hsize 3cm\parindent=0pt\eightpoints\baselineskip=9pt%
      \rightskip=4mm plus 4mm\raggedright#1\hss%
      }%
    \hbox{\kern-4cm\smash{\box1}\hss\par}%
    }%
  \hfuzz=0pt
  }
\def\note#1{\relax}

\def\anote#1#2#3{\smash{\kern#1in{\raise#2in\hbox{#3}}}%
  \nointerlineskip}     

\newcount\equanumber
\equanumber=0
\newcount\sectionnumber
\sectionnumber=0
\newcount\subsectionnumber
\subsectionnumber=0
\newcount\snumber  
\snumber=0

\def\section#1{%
  \subsectionnumber=0%
  \snumber=0%
  \equanumber=0%
  \advance\sectionnumber by 1%
  \noindent{\bf \the\sectionnumber .~#1.~}%
}%
\def\subsection#1{%
  \advance\subsectionnumber by 1%
  \snumber=0%
  \equanumber=0%
  \noindent{\bf \the\sectionnumber .\the\subsectionnumber .~#1.~}%
}%
\def\prevs{\the\sectionnumber .\the\subsectionnumber .\the\snumber }

\long\def\Corollary#1{%
  \global\advance\snumber by 1%
  \bigskip
  \noindent{\bf Corollary~\prevs .}%
  \quad{\it#1}%
}%
\long\def\Lemma#1{%
  \global\advance\snumber by 1%
  \bigskip
  \noindent{\bf Lemma~\prevs .}%
  \quad{\it#1}%
}%
\def\Proof{\noindent{\bf Proof.~}}
\long\def\Proposition#1{%
  \advance\snumber by 1%
  \bigskip
  \noindent{\bf Proposition~\prevs .}%
  \quad{\it#1}%
}%
\long\def\Theorem#1{%
  \advance\snumber by 1%
  \bigskip
  \noindent{\bf Theorem~\prevs .}%
  \quad{\it#1}%
}%
\long\def\Statement#1{%
  \advance\snumber by 1%
  \bigskip
  \noindent{\bf Statement~\prevs .}%
  \quad{\it#1}%
}%
\def\ifundefined#1{\expandafter\ifx\csname#1\endcsname\relax}
\def\labeldef#1{\global\expandafter\edef\csname#1\endcsname{\prevs}}
\def\labelref#1{\expandafter\csname#1\endcsname}
\def\label#1{\ifundefined{#1}\labeldef{#1}\note{$<$#1$>$}\else\labelref{#1}\fi}

\def\preveq{(\the\sectionnumber .\the\subsectionnumber .\the\equanumber)}
\def\neq{\global\advance\equanumber by 1\eqno{\preveq}}

\def\ifundefined#1{\expandafter\ifx\csname#1\endcsname\relax}

\def\equadef#1{\global\advance\equanumber by 1%
  \global\expandafter\edef\csname#1\endcsname{\preveq}%
  \preveq}

\def\equaref#1{\expandafter\csname#1\endcsname}

\def\equa#1{%
  \ifundefined{#1}%
    \equadef{#1}%
  \else\equaref{#1}\fi}

\font\eightrm=cmr8%
\font\sixrm=cmr6%

\font\eightsl=cmsl8%

\font\eightbf=cmb8%

\font\eighti=cmmi8%
\font\sixi=cmmi6%

\font\eightsy=cmsy8%
\font\sixsy=cmsy6%

\font\eightex=cmex8%
\font\sixex=cmex6%
\font\fiveex=cmex5%

\font\eightit=cmti8%

\font\eighttt=cmtt8%

\font\tenbb=msbm10%
\font\eightbb=msbm8%
\font\sevenbb=msbm7%
\font\sixbb=msbm6%
\font\fivebb=msbm5%
\newfam\bbfam  \textfont\bbfam=\tenbb  \scriptfont\bbfam=\sevenbb  \scriptscriptfont\bbfam=\fivebb%

\font\tenbbm=bbm10

\font\tencmssi=cmssi10%
\font\sevencmssi=cmssi7%
\font\fivecmssi=cmssi5%
\newfam\ssfam  \textfont\ssfam=\tencmssi  \scriptfont\ssfam=\sevencmssi  \scriptscriptfont\ssfam=\fivecmssi%
\def\ssi{\fam\ssfam\tencmssi}%

\font\tenfrak=cmfrak10%
\font\eightfrak=cmfrak8%
\font\sevenfrak=cmfrak7%
\font\sixfrak=cmfrak6%
\font\fivefrak=cmfrak5%
\newfam\frakfam  \textfont\frakfam=\tenfrak  \scriptfont\frakfam=\sevenfrak  \scriptscriptfont\frakfam=\fivefrak%
\def\frak{\fam\frakfam\tenfrak}%

\font\tenmsam=msam10%
\font\eightmsam=msam8%
\font\sevenmsam=msam7%
\font\sixmsam=msam6%
\font\fivemsam=msam5%

\def\bb{\fam\bbfam\tenbb}%

\def\hexdigit#1{\ifnum#1<10 \number#1\else%
  \ifnum#1=10 A\else\ifnum#1=11 B\else\ifnum#1=12 C\else%
  \ifnum#1=13 D\else\ifnum#1=14 E\else\ifnum#1=15 F\fi%
  \fi\fi\fi\fi\fi\fi}
\newfam\msamfam  \textfont\msamfam=\tenmsam  \scriptfont\msamfam=\sevenmsam  \scriptscriptfont\msamfam=\fivemsam%
\def\msam{\msamfam\tenmsam}%
\mathchardef\leq"3\hexdigit\msamfam 36%
\mathchardef\geq"3\hexdigit\msamfam 3E%

\font\tentt=cmtt11%
\font\seventt=cmtt9%
\textfont\ttfam=\tentt
\scriptfont7=\seventt%
\def\tt{\fam\ttfam\tentt}%

\def\eightpoints{%
\def\rm{\fam0\eightrm}%
\textfont0=\eightrm   \scriptfont0=\sixrm   \scriptscriptfont0=\fiverm%
\textfont1=\eighti    \scriptfont1=\sixi    \scriptscriptfont1=\fivei%
\textfont2=\eightsy   \scriptfont2=\sixsy   \scriptscriptfont2=\fivesy%
\textfont3=\eightex   \scriptfont3=\sixex   \scriptscriptfont3=\fiveex%
\textfont\itfam=\eightit  \def\it{\fam\itfam\eightit}%
\textfont\slfam=\eightsl  \def\sl{\fam\slfam\eightsl}%
\textfont\ttfam=\eighttt  \def\tt{\fam\ttfam\eighttt}%
\textfont\bffam=\eightbf  \def\bf{\fam\bffam\eightbf}%

\textfont\frakfam=\eightfrak  \scriptfont\frakfam=\sixfrak \scriptscriptfont\frakfam=\fivefrak  \def\frak{\fam\frakfam\eightfrak}%
\textfont\bbfam=\eightbb      \scriptfont\bbfam=\sixbb     \scriptscriptfont\bbfam=\fivebb      \def\bb{\fam\bbfam\eightbb}%
\textfont\msamfam=\eightmsam  \scriptfont\msamfam=\sixmsam \scriptscriptfont\msamfam=\fivemsam  \def\msam{\msamfam\eightmsam}

\rm%
}

\def\poorBold#1{\setbox1=\hbox{#1}\wd1=0pt\copy1\hskip.25pt\box1\hskip .25pt#1}

\mathchardef\lsim"3\hexdigit\msamfam 2E%
\mathchardef\gsim"3\hexdigit\msamfam 26%

\def\d{\,{\rm d}}

\def\ds{\displaystyle}

\long\def\DoNotPrint#1{\relax}
\def\E{{\rm E}}
\def\eqd{\buildrel {\rm d}\over =}

\def\existsnl{\exists n\in(\Lambda,\Lambda^{1+\epsilon})\,:\,}

\def\fixedref#1{#1\note{fixedref$\{$#1$\}$}}

\def\g#1{g_{[0,#1)}}
\def\Id{{\rm Id}}
\def\kl{k(\Lambda)}

\def\lima{\lim_{a\to 0}}

\def\liminft{\liminf_{t\to\infty}}

\def\liml{\lim_{\Lambda\to\infty}}
\def\limn{\lim_{n\to\infty}}

\def\limsupL{\limsup_{\Lambda\to\infty}}

\def\limsupt{\limsup_{t\to\infty}}
\def\limt{\lim_{t\to\infty}}
\def\limT{\lim_{T\to\infty}}

\def\ml{m_\Lambda}

\def\oF{\overline F{}}

\def\prob#1{\hbox{\rm P}\{\,#1\,\}}
\def\proB#1{\hbox{\rm P}\bigl\{\,#1\,\bigr\}}
\def\Prob#1{\hbox{\rm P}\Bigl\{\,#1\,\Bigr\}}
\def\qed{{\vrule height .9ex width .8ex depth -.1ex}}
\def\ROne{{\cal R}_{1,n,N}}

\def\ss{\scriptstyle}

\def\Var{{\rm Var}}

\def\tough{Barbe and McCormick (2010b)}
\def\gFLevy{Barbe and McCormick (2010a)}

\def\boc{\note{{\bf BoC}\hskip-11pt\setbox1=\hbox{$\Bigg\downarrow$}%
         \dp1=0pt\ht1=0pt\ht1=0pt\leavevmode\raise -20pt\box1}}
\def\eoc{\note{{\bf EoC}\hskip-11pt\setbox1=\hbox{$\Bigg\uparrow$}%
         \dp1=0pt\ht1=0pt\ht1=0pt\leavevmode\raise 20pt\box1}}

\def\One{\hbox{\tenbbm 1}}

\def\MM{{\bb M}\kern .4pt}
\def\NN{{\bb N}\kern .5pt}
\def\RR{{\bb R}}

\def\UU{{\bb U}}

\def\M{{\ssi M}}


\pageno=1
\centerline{\bf RUIN PROBABILITIES IN TOUGH TIMES}
\vskip 8pt
\centerline{\bf\poorBold{$\widetilde{\hbox to 1cm{\hfill}}$}}
{\eightpoints
\centerline{\bf Part 2}
\centerline{\bf HEAVY-TRAFFIC APPROXIMATION FOR FRACTIONALLY}
\centerline{\bf DIFFERENTIATED RANDOM WALKS IN THE DOMAIN OF ATTRACTION}
\centerline{\bf OF A NONGAUSSIAN STABLE DISTRIBUTION}
}
\bigskip

\centerline{Ph.\ Barbe$^{(1)}$ and W.P.\ McCormick$^{(2)}$}
\centerline{${}^{(1)}$CNRS {\sevenrm(UMR {\eightrm 8088})}, ${}^{(2)}$University of Georgia}

{\narrower
\baselineskip=9pt\parindent=0pt\eightpoints

\bigskip

{\bf Abstract.} Motivated by applications to insurance mathematics, we prove
some heavy-traffic limit theorems for processes which encompass the 
fractionally differentiated random walk as well as some FARIMA processes, 
when the innovations are in the domain of attraction of a nonGaussian stable 
distribution.

\bigskip


\bigskip

\hfuzz=3pt
\noindent{\bf Keywords:} heavy traffic, ruin probability, fractional random
walk, FARIMA process, Poisson process.

\hfuzz=0pt
}

\bigskip\bigskip


\def\preveq{(\the\sectionnumber.\the\equanumber)}
\def\prevs{\the\sectionnumber.\the\snumber }

\section{Introduction and main result}%
The purpose of this paper is to complete our study of ruin probability for
some nonstationary processes having long range dependence and innovations in 
the domain of attraction of a nonGaussian stable distribution, when the
premium can hardly cover the claims. The overall motivations for this study
are described in the first part of this work (Barbe and McCormick, 2010b).

To recall the setting, we consider a function $g$ analytic on $(-1,1)$,
$$
  g(x)=\sum_{i\geq 0} g_i x^i \, .
$$
Given a distribution function $F$, we can define a $(g,F)$-process $(S_n)$ as 
follows.
Let $(X_i)_{i\geq 1}$ be a sequence of random variables, independent, having
common distribution function $F$. We set $X_i=0$ if $i$ is nonpositive. We
define the backward shift operator $B$ by $BX_i=X_{i-1}$ and set
$$
  S_n = g(B)X_n = \sum_{0\leq i<n} g_i X_{n-i} \, .
$$
Important examples of such processes include 

{
\setbox1=\hbox{--}
\setbox1=\hbox{--~}
\parindent=-\wd1
\leftskip\wd1

\smallskip

--~random walks: $g(x)=(1-x)^{-1}$, 

\smallskip

--~ARMA processes: $g$ a rational function
with poles outside the complex unit disk, and

\smallskip

--~FARIMA processes: $g(x)=(1-x)^{-\gamma} R(x)$ where $\gamma$ is positive 
and $R$ is a rational function not vanishing at $1$ and with poles outside 
the complex unit disk.\par
}

\medskip

\noindent
In particular, if $g(x)=(1-x)^{-\gamma}$, then 
$$
  S_n=(1-B)^{1-\gamma}(1-B)^{-1}X_n
$$ 
is the random walk $(1-B)^{-1}X_n$ differentiated $1-\gamma$ times.

Following the first part of this work, we assume that the sequence 
\setbox1=\vbox{\hsize=3in\noindent $(g_n)$ is ultimately positive and regularly varying 
of negative index $\gamma-1$ with $\gamma$ in $(0,1)$.}
$$
  \box1
  \eqno{\equa{hypRVg}}
$$
This forces the function $g$ to diverge to $+\infty$ as its argument tends to
$1$.
Part 1 of this work concentrates on the case where $\gamma$ is greater 
than $1$, forcing $(g_n)$ to diverge to infinity. In contrast, in this
current paper, we assume that $\gamma$ is less than $1$, forcing $(g_n)$
to converge to $0$.

For any positive $t$ define the partial sum
$$
  \g t=\sum_{0\leq i<t} g_i \, .
$$
Under \hypRVg, the sequence $(\g n)$ diverges to infinity. Moreover, 
under \hypRVg, $(g_n)$ is asymptotically 
equivalent to a monotone sequence, and Karamata's theorem for power series 
(see Bingham, Goldie and Teugels, 1989, Corollary 1.7.3) asserts that
$$
  \g n\sim {ng_n\over\gamma}\sim {g(1-1/n)\over\Gamma(1+\gamma)}
  \eqno{\equa{gEquiv}}
$$
as $n$ tends to infinity. In particular, writing $\Id$ for the identity
function on the real line, $g(1-1/\Id)$ is regularly varying of index $\gamma$
at infinity.

If $\E X_1$ is negative, $\E S_n=\g n \E X_1$ diverges toward minus infinity 
as $n$ tends to infinity. It is then conceivable that $\sup_n S_n$ might be
almost surely finite, and, if this is the case, our heavy traffic 
approximation describes the limiting behavior of this all time supremum when
the expectation of the innovations tends to $0$. Writing $S_n=(S_n-\g n \E X_1)
+\E X_1 \g n$, an alternative viewpoint is to consider that
$$
  \hbox{$F$ is centered}
  \eqno{\equa{hypCentered}}
$$
and seek the limiting behavior of $\sup_{n\geq 0} (X_n-a\g n)$ when $a$ tends
to $0$. As indicated in the first part of this work, this problem has bearing
on calculations of ruin probabilities in insurance risk when premiums can 
barely keep up with
claims, in queueing theory, and for moving boundary crossing probabilities as
well.

Throughout the paper we will use $c$ for a generic constant whose value
may change from place to place.

As in the first part, we assume that
\setbox1=\vbox{\hsize=210pt\noindent%
  $F$ belongs to the domain of attraction of a stable
  distribution with index $\alpha$ in $(1,2)$.}
$$
  \box1
  \eqno{\equa{hypStable}}
$$

Assumption \hypStable\ implies that $F$ is tail balanced in the following 
sense. Writing $F_*$ for the distribution of $|X_1|$ and $\M_{-1}F$ for that of
$-X_1$, there exist $p$ and $q$ both in $[\,0,1\,]$, such that
$\oF\sim p \oF_*$ and $\overline{\M_{-1}F}\sim q\oF_*$ at infinity. These
asymptotic relations imply $p+q=1$. Under \hypStable, $\oF_*$ is regularly 
varying of index $-\alpha$; and so is $\oF$, $\overline{\M_{-1}F}$,
if $p$, $q$, does not vanish respectively. As in the first part of this work,
we will assume throughout this paper that $p$ does not vanish. If $q$ vanishes
we will also assume that the lower tail of the distribution function $F$ decays
slightly faster than the upper one in the sense that for some constant $c$
$$
  \overline{\M_{-1}F}(t)\leq c \oF(t\log t)\log t\qquad\hbox{ultimately.}
  \eqno{\equa{tailDomin}}
$$
Assumption \hypStable\ gives rise to a L\'evy measure $\nu$ whose density
with respect to the Lebesgue measure $\lambda$ is
$$
  {\d\nu\over\d\lambda}(x)
  = p\alpha x^{-\alpha-1} \One_{(0,\infty)}(x) 
    + q\alpha (-x)^{-\alpha-1}\One_{(-\infty,0)}(x) \, .
$$

The function
$$
  k(t)
  ={\g t\over F_*^\leftarrow(1-1/t)}
$$
will play a role in our following main result ---~note that the meaning of
the notation $k$ in this paper is different than that in part I, but will play
an analogous role. Given \gEquiv, it is asymptotically equivalent 
to $g(1-1/t)/\bigl(\Gamma(1+\gamma)F_*^\leftarrow(1-1/t)\bigr)$ as $t$ tends
to infinity.

The first assertion of our main theorem asserts that for the heavy traffic
to make sense we should have $\alpha\gamma\geq 1$. 

\Theorem{\label{mainTh}
  Assume that \hypRVg\ holds for some positive $\gamma$ less than $1$, and 
  that \hypCentered\ and \hypStable\ hold. If $q$ vanishes, assume furthermore
  that \tailDomin\ holds. Then

  \medskip

  \noindent (i) if $\limt k(t)<\infty$, in particular
  if $\alpha\gamma<1$, then for any positive $a$,
  $$
    \sup_{n\geq 0}S_n-a\g n=+\infty
  $$
  in probability.

  \noindent (ii) if $\alpha\gamma>1$ then, writing 
  $\sum_{i\geq 1}\delta_{(t_i,x_i)}$ for a Poisson process with mean intensity
  $\lambda\otimes\nu$, the distribution of 
  $$
    {\sup_{n\geq 0} (S_n-a\g n)\over 
       F_*^\leftarrow\Bigl(1-{\ds 1\over\ds k^\leftarrow(1/a)}\Bigr) } 
  $$
  converges to that of
  $
    \sup_{\ss i\geq 1\atop\ss j\geq 0} g_j x_i-t_i^\gamma\, .
  $
  The latter random variable is almost surely finite.
}

\bigskip

Note that Theorem \mainTh\ leaves open the boundary case where $\alpha\gamma=1$
and $\limt k(t)=+\infty$. It is conceivable
that the conclusion of (ii) remains, but we do not know how to prove it.
Theorem \mainTh\ also leaves open the seemingly less interesting situation 
where $k$ oscillates in such a way that $\liminft k(t)=0$ 
and $\limsupt k(t)=+\infty$.

The following elementary example gives a concrete idea on the rate involved 
in Theorem \mainTh. Using the first part of this work, it is straightforward
to extend this example to FARIMA models.

\medskip

\noindent{\it Example.} Consider $g_i=i^{\gamma-1}$ with $\gamma<1$. Assume
that $\oF_*(t)\sim c t^{-\alpha}$ as $t$ tends to infinity, and 
that $\alpha\gamma>1$. This implies $\oF_*^\leftarrow(1-1/t)\sim 
(ct)^{1/\alpha}$ as $t$ tends to infinity. From \gEquiv\ we
deduce that $\g i\sim i^\gamma/\gamma$. Thus, $k(t)\sim t^{\gamma-1/\alpha}/
(\gamma c^{1/\alpha})$ and $k^\leftarrow(1/a)
\sim (\gamma c^{1/\alpha}/a)^{\alpha/(\alpha\gamma-1)}$ as $a$ tends to $0$.
We then obtain $F_*^\leftarrow\bigl(1-1/k^\leftarrow(1/a)\bigr)\sim 
(\gamma c^\gamma/a)^{1/(\alpha\gamma-1)}$. Assertion (ii) of Theorem \mainTh\ 
implies that
$\sup_{n\geq 0} (S_n-a \g n)$ grows like $1/a^{1/(\alpha\gamma-1)}$ as
$a$ tends to $0$.

\bigskip


\section{Proof of Theorem \mainTh}
Let $\Pi=\sum_i\delta_{(t_i,x_i)}$ be a Poisson random measure with mean
intensity $\lambda\otimes\nu$. For any positive integer $i$, set
$$
  c_{n,i}=\g i \int_{F^\leftarrow(1/n)}^{F^\leftarrow(1-1/n)} x \d F(x) \, .
$$
Since $F$ is centered, 
$$
  \int_{F^\leftarrow(1/n)}^{F^\leftarrow(1-1/n)} x \d F(x)
  = O\bigl(F_*^\leftarrow(1-1/n)/n\bigr)
$$
as $n$ tends to infinity. Thus, since $\gamma$ is less than $1$, for 
any positive $M$,
$$
  \limn \max_{1\leq i\leq Mn} c_{n,i}/F_*^\leftarrow(1-1/n)
  =0 \, .
  \eqno{\equa{cinZero}}
$$
Consider the random measure 
$$
  N_n=\sum_{i\geq 1} \delta_{(i/n,S_i/F_*^\leftarrow(1-1/n))}
$$
in the space of all measures on $[\,0,\infty)\times (\RR\setminus\{\,0\,\})$
endowed with the topology of vague convergence.
Theorem 5.3 in Barbe and McCormick (2010a)
asserts that the distribution of the random measure 
$$
  \sum_{i\geq 1}\delta_{(i/n,(S_i-c_{n,i})/F_*^\leftarrow(1-1/n))}
$$ 
converges weakly$*$ to that of 
$N=\sum_{i\geq 1} \sum_{j\geq 0} \delta_{(t_i,g_j x_i)}$. Since \cinZero\ 
holds, this implies that the distribution of $N_n$ converges to that of $N$
as well.

Define $\Lambda=\Lambda(1/a)$ by
$$
  k(\Lambda)\sim 1/a
$$
as $a$ tends to $0$, that is $\Lambda\sim k^\leftarrow$. Let $T$ be 
a positive real number. We obtain
$$\eqalign{
  {S_i-a\g i\over F_*^\leftarrow(1-1/\Lambda)}
  &{}={S_i-a\g\Lambda {\ds\g i\over\ds\g\Lambda}\over F_*^\leftarrow(1-1/\Lambda)}
   \cr
  &{}={S_i\over F_*^\leftarrow(1-1/\Lambda)} - {\g i\over\g\Lambda}
    \bigl(1+o(1)\bigr) \cr
  }
$$
where the $o(1)$ term is uniform in $i$ between $0$ and $\Lambda T$ and as
$a$ tends to $0$.

Let $\epsilon$ be a positive real number. By the Skorokhod-Dudley-Wichura
representation theorem, we can construct a version of $N$ and, for each $n$,
a version of $N_n$ such that this version of $N_n$ converges almost surely
to $N$ as point measures on $[\,0,T\,]\times (\RR\setminus\{\,0\,\})$.
We consider these versions until subsection \fixedref{2.2}, even though
we use the same notation as the original processes. For these versions,
$$
  \sup_{i\geq 0} {S_i\One\{\, S_i/F_*^\leftarrow(1-1/\Lambda)>\epsilon\,\}
                  -a\g i\over F_*^\leftarrow(1-1/\Lambda)}
  \One\{\, 0\leq i\leq \Lambda T \,\}
  \eqno{\equa{ThBEqA}}
$$
converges almost surely to
$$
  \sup_{i\geq 1} \sup_{j\geq 0} 
  (g_jx_i\One\{\, g_jx_i>\epsilon\,\}-t_i^\gamma)
  \One\{\, 0\leq t_i\leq T\,\} \, .
$$
Since $\epsilon$ is arbitrary and 
$\sup_{0\leq i\leq \Lambda T} (S_i-a\g i)/F_*^\leftarrow(1-1/\Lambda)$
is within $\epsilon$ of \ThBEqA, this implies
$$
  \lima \sup_{0\leq i\leq \Lambda T} {S_i-a\g i\over F_*^\leftarrow(1-1/\Lambda)}
  = \sup_{i\geq 1} \sup_{j\geq 0} (g_jx_i-t_i^\gamma)
  \One\{\, 0\leq t_i\leq T\,\}  \, .
  \eqno{\equa{ThBEqB}}
$$

With these preliminaries, we can prove both assertions of Theorem \mainTh\ in
the next two subsections.

\bigskip

\def\preveq{(\the\sectionnumber .\the\subsectionnumber .\the\equanumber)}
\def\prevs{\the\sectionnumber .\the\subsectionnumber .\the\snumber }

\subsection{Proof of assertion (i)} If $x$ is a real number, we write $x_+$ for
$x\vee 0$. For any nonnegative integer $p$, define
$$
  M_p=\max_{i:p\leq t_i<p+1} \sup_{j\geq 0}(g_j x_i)_+ \, .
$$
Let $T$ be a positive integer. Since our version of $N_n$ converges almost
surely to $N$,
$$\displaylines{\quad
  \max_{0\leq i<nT} (S_i)_+/\g i
  \hfill\cr\hfill
  \eqalign{
  {}={}&\max_{0\leq p<T} \max_{np\leq i<n(p+1)} 
        { (S_i)_+\over\ds F_*^\leftarrow(1-1/n)}
          {\ds F_*^\leftarrow(1-1/n)
        \over\ds F_*^\leftarrow(1-1/i)}
        {1\over k(i)} \cr
  {}\geq{}& \max_{1\leq p<T} \bigl(M_p+o(1)\bigr) 
            {1\over (p+1)^{1/\alpha}}
            {1+o(1)\over\limsupt |k(t)|} \cr}
  \quad\cr}
$$
almost surely. It then suffices to show that 
$$
  \limT \max_{1\leq p\leq T}M_p/p^{1/\alpha} =+\infty
$$ 
in probability. Since $\gamma$ is less than $1$, the sequence
$(g_j)_{j\geq 0}$ has a largest term, $g_\infty$, which is positive 
under \hypRVg. Let $(\omega_i)_{i\geq 1}$ be a sequence of independent random
variables having an exponential distribution with mean $1$. The discussion
following Lemma 6.1 in \gFLevy, or the calculation between \fixedref{(2.2.11)}
and \fixedref{(2.2.12)} in this paper, show that $(M_p)_{p\geq 1}$ has
the same distribution as $(g_\infty \omega_p^{-1/\alpha})_{p\geq 1}$. Thus, it 
suffices to show that $\min_{p\geq 1} p\omega_p=0$ in probability. This
follows from the equality
$$
  \proB{ \min_{1\leq p\leq T} p\omega_p >\epsilon}
  =\prod_{1\leq p\leq T} e^{-\epsilon/p}
$$
and the divergence of the series $\sum_{p\geq 1} 1/p$.

\bigskip

\subsection{Proof of assertion (ii)}
Given our preliminary remarks, and in particular \ThBEqB,
it suffices to prove that
$$
  \limT\limsupL \prob{\exists n>\Lambda T\,:\, S_n>a\g n} = 0 \, .
$$
Arguing as in the beginning of section 3 of \tough,
this is equivalent to showing that
$$
  \limT\limsupL\prob{\exists n>\Lambda \,:\, S_n>T\g n/\kl} = 0 \, .
  \eqno{\equa{eqKey}}
$$

This will be achieved in mostly three steps, and a fourth one to handle the 
part of the proof dealing with the lower tail of the distribution.

\noindent{\it Step 1.}
Let $(a_n)$ and $(b_n)$ be two sequences diverging respectively to $-\infty$
and $+\infty$. We assume that
$$
  \limn b_n/(-a_n) \hbox{ is positive or infinite.}
$$
We set
$$
  \sigma_n^2=\Var\bigl(X_1\One_{[a_n,b_n]}(X_1)\bigr)
$$
and, for any positive integer $i$ at most $n$,
$$
  Z_{i,n}={X_i\One_{[a_n,b_n]}(X_i)-\E X_i\One_{[a_n,b_n]}(X_i)\over \sigma_n}
  \, . 
$$
Up to scaling, the part of $S_n$ made by the `middle' innovations is 
$$
  M_n=\sigma_n\sum_{0\leq i <n} g_i Z_{n-i,n} \, .
$$
As in \tough, we construct $(b_n)$ to be regularly varying, of
the form $F^\leftarrow(1-m_n/n)$ with $1-m_n/n$ in the range of $F$ and $(m_n)$
being regularly varying of index $\beta$.

\Proposition{\label{middleNeglect}
  For any positive $\beta$ less than $1$ and any positive $T$,
  $$
    \liml\prob{\exists n\geq\Lambda \,:\, M_n>T\g n/\kl} = 0 \, .
  $$
}

\Proof As in the proof of Proposition 3.1.1 in \tough,
$$
  \sigma_n\sim c F^\leftarrow(1-m_n/n)\sqrt{m_n/n}
$$
as $n$ tends to infinity. Moreover, for any positive integer $r$ such 
that $\sum |g_i|^r$ converges,
$$
  \Bigl| \E\Bigl({M_n\over \sigma_n\sqrt n}\Bigr)^r\Bigr|
  \leq {c_r\over n} \, .
$$
Markov's inequality yields
$$\eqalign{
  \Prob{ M_n>T\g n/\kl}
  &{}\leq \Bigl({\sigma_n \sqrt{n}\kl\over T\g n}\Bigr)^r {c_r\over n} \cr
  &{}\sim {c\over T^r} \Bigl({F^\leftarrow(1-m_n/n)\sqrt{m_n}\kl\over \g n}
    \Bigr)^r {1\over n} \cr
  &{}\sim {c\over T^r} \kl^r \Bigl( {F^\leftarrow(1-m_n/n)\sqrt{m_n}\over \g n}
    \Bigr)^r {1\over n}\, .\cr
  }
$$
This asymptotic equivalent is of the form $\kl^r$ times a function of $n$ 
which is regularly varying of index
$$
  r\Bigl( {1-\beta\over\alpha}+{\beta\over 2}-\gamma\Bigr) -1
  = r\Bigl( {1\over\alpha}-\gamma+\beta\Bigl({1\over 2}-{1\over\alpha}\Bigr)
  \Bigr) - 1 \, .
$$
This index is less than $-1$. Therefore, by Bonferroni's
inequality and Karamata's theorem,
$$\displaylines{\qquad
  \prob{\exists n>\Lambda \,:\, M_n>T \g n/\kl}
  \hfill\cr\hfill
  \leq {c\over T^r} \kl^r
  \Bigl({F^\leftarrow(1-\ml/\Lambda)\sqrt{\ml}\over\g\Lambda}\Bigr)^r \, .
  \qquad\cr}
$$
This bound is regularly varying in $\Lambda$ of negative 
index $r\beta\bigl((1/2)-(1/\alpha)\bigr)$ and tends to $0$ as $\Lambda$ 
tends to infinity.\hfill\qed

\bigskip

\noindent{\it Step 2.} We consider the part of $S_n$ made by the extreme
innovations,
$$
  T_n^+=\sum_{0\leq i <n} g_i X_{n-i}\One_{[b_n,\infty)}(X_{n-i}) \, .
$$
The purpose of this step is to show that in our problem we can ignore the
contribution of $T_n^+$ in the range of $n$ exceeding $\Lambda^{1+\epsilon}$.
The following lemma will be instrumental; it is stronger than what we need
in this step, but this strength will turn useful in the next step.

\Lemma{\label{expectNeglect}
  If $\beta<(1-\gamma)/(1-1/\alpha)$ then whenever $T$ is large enough, $n$
  is at least $\Lambda$ and $\Lambda$ is large enough, $\E T_n^+\leq
  T\g n/\kl$.
}

\bigskip

\Proof
Lemma 3.2.2 in \tough\ implies
$$
  \E T_n^+ \sim \g n {\alpha\over\alpha-1} m_n^{1-1/\alpha} 
  {F^\leftarrow(1-1/n)\over n}\, .
$$
Substituting $T$ by a multiple of it, it suffices to prove that
$$
  {\g n\over n} m_n^{1-1/\alpha} < T {k(n)\over\kl} \, .
$$
Since $k$ is regularly varying of positive index, it suffices to show that
$$
  {\g n \over n} m_n^{1-1/\alpha} < T \, .
$$
This holds because the left hand side is regularly varying of index
$$
  \gamma-1+\beta(1-1/\alpha)
$$
which is negative.\hfill\qed

\bigskip

The main result of this step 2 is the following.

\Proposition{\label{stepTwo}
  Let $\epsilon$ be a positive real number. If 
  $$
    \beta<{1\over 2\gamma}\Bigl(\gamma-{1\over \alpha}\Bigr) 
    {\epsilon\over 1+\epsilon} \, ,
  $$
  then
  $$
    \limT\limsupL\prob{\exists n>\Lambda^{1+\epsilon} \, : \, |T_n^+-\E T_n^+|
    > T\g n/\kl} = 0 \, .
  $$
}

\Proof Let $(U_i)_{i\geq 1}$ be a sequence of independent random variables
uniform over $[\,0,1\,]$. Let $\UU_n$ be the empirical distribution function
of $(U_i)_{1\leq i\leq n}$. We write $(U_{i,n})$ for the order statistics of
$(U_i)_{1\leq i\leq n}$. Without any loss of generality, we assume that
$X_i=F^\leftarrow(1-U_i)$, so that, for any $n$ large enough
$$\eqalign{
  T_n^+
  &{}=\sum_{0\leq i<n} g_i F^\leftarrow(1-U_{n-i})
    \One\{\, U_{n-i}\leq m_n/n\,\} \cr
  &{}\leq c F^\leftarrow(1-U_{1,n})\g{n\UU_n(m_n/n)} \, .\cr
  }
$$

Let $(\xi_n)_{n\geq 1}$ be a slowly varying nondecreasing sequence such that
$\sum_{n\geq 1} 1/n\xi_n$ converges.
From Kiefer's (1972) Theorem 1 we deduce that $U_{1,n}\geq 1/\xi_n$ almost
surely for $n$ large enough, while Shorack and Wellner's (1978) Theorem 2 
implies $\UU_n
\leq \xi_n \Id$. Then, using Potter's bound for $F^\leftarrow(1-1/\Id)$
and using that $(\g n)$ is regularly varying, we obtain
$$
  T_n^+\leq c F^\leftarrow(1-1/n)\xi_n^{2/\alpha} \g{m_n} \, .
$$
Therefore, the inequality $T_n^+>T \g n/\kl$ implies, almost surely for $n$
large enough,
$$
  \kl > cT {k(n)\over \g{m_n}\xi_n^{2/\alpha}} \, .
  \eqno{\equa{stepTwoA}}
$$
In this inequality, the right hand side is regularly varying of index 
$\gamma-(1/\alpha)-\beta\gamma$, which is positive provided 
$\beta<\bigl(\gamma-(1/\alpha)\bigr)/\gamma$. In this range of $\beta$ and in 
the range of $n$ at least $\Lambda^{1+\epsilon}$, the right hand side of
\stepTwoA\ is at least a constant times its value at $\Lambda^{1+\epsilon}$.
This lower bound, as a function of $\Lambda$, is regularly varying
and for \stepTwoA\ to hold we must have, comparing the index of regular 
variation,
$$
  \gamma-{1\over\alpha}
  \geq (1+\epsilon)\Bigl(\gamma-{1\over\alpha}-\beta\gamma\Bigr) \, .
$$
This does not hold under the assumption of the lemma, and therefore \stepTwoA\
does not occur. So $T_n^+\leq T\g n/\kl$ almost surely in the range 
$n\geq \Lambda^{1+\epsilon}$ and for $\Lambda$ large enough.

Lemma \expectNeglect\ implies that
$$
  \E T_n^+ \leq T \g n /\kl \, ,
$$
whenever $n$ exceeds $\Lambda^{1+\epsilon}$ and $\Lambda$ is large enough.
This proves the proposition.\hfill\qed

\bigskip

\noindent{\it Step 3.} Given Lemma \expectNeglect, our goal is now to show that
$$
  \limT\limsupL\Prob{\existsnl T_n^+> T\g n/\kl} = 0 \, .
  \eqno{\equa{stepThreeA}}
$$
For this, we approximate $T_n^+$ by a simpler quantity.

It is convenient to write $N$ for $\Lambda^{1+\epsilon}$.
Let $(U_i)_{i\geq 1}$ be a sequence of independent random variables uniformly
distributed on $[\,0,1\,]$. Let $\tau$ be the random permutation of $\{\,1 ,2,
\ldots,N\,\}$ such that $U_{\tau(i)}=U_{i,N}$. Without any loss of generality
we can assume that $X_i=F^\leftarrow(1-U_i)$. For $n$ in 
$(\Lambda,\Lambda^{1+\epsilon})$, we then have
$$\eqalign{
  T_n^+
  &{}=\sum_{1\leq i\leq n} g_{n-i} X_i\One\{\, X_i >b_n\,\} \cr
  &{}=\sum_{1\leq i\leq N} g_{n-\tau(i)} F^\leftarrow(1-U_{i,N}) 
    \One\{\, U_{i,N}\leq m_n/n\,\} \One\{\, \tau(i)\leq n\,\} \, . \cr
  }
$$
Let $(V_i)_{i\geq 1}$ be a sequence of independent random variables uniformly
distributed on $[\,0,1\,]$, independent of $(U_i)$. Let $G_N$ be the 
empirical distribution function of $(V_i)_{1\leq i\leq N}$. Without any loss
of generality we can assume that $\tau(i)=NG_N(V_i)$. Then, setting
$$
  \ROne=\Bigl\{\, i\,:\, U_{i,N}\leq {m_n\over n} \, ;\, 
                  G_N(V_i)\leq {n\over N} \, \Bigr\} \, ,
$$
we have
$$
  T_n^+ = \sum_{i\in \ROne} g_{n-NG_N(V_i)} F^\leftarrow(1-U_{i,N}) \, .
$$
Let $i^*=i^*_{n,N}$ be in $\ROne$ such that $n-NG_N(V_i)$ is minimum when
$i$ is
$i^*$. Such $i^*$ exists and is well defined because $\ROne$ is finite and
for $i$ in $\ROne$ the differences $n-NG_N(V_i)$ are nonnegative and assume,
almost surely, different values for different $i$. Set
$$
  T_{1,n,N}^+=g_{n-NG_N(V_{i^*})} F^\leftarrow(1-U_{i^*,N}) \, .
$$
Our next lemma shows that we can approximate $T_n^+$ by $T_{1,n,N}^+$ in 
order to prove \stepThreeA.

\Lemma{\label{approxOne}
  For any $\epsilon$ and $\beta$ small enough, for any positive $T$,
  $$
    \liml \prob{\existsnl |T_n^+-T_{1,n,N}^+|>T\g n/\kl} = 0 \, .
  $$
}

\Proof Writing $\sharp \ROne$ for the cardinality of $\ROne$, we have
$$\displaylines{
  |T_n^+-T_{1,n,N}^+|
  \hfill\cr\noalign{\vskip 3pt}\hfill
    \eqalign{
      {}={}&    \sum_{i\in\ROne\setminus\{i^*\}} g_{n-NG_N(V_i)} 
                F^\leftarrow(1-U_{i,N}) \cr
      {}\leq{}& \sharp\ROne \max_{i\in\ROne\setminus\{i^*\}} g_{n-NG_N(V_i)}
                \max_{i\in\ROne} F^\leftarrow(1-U_{i,N}) \, . 
                \ \equa{approxOneB}\cr}
  \cr}
$$
Let $\eta$ be a positive real number less than $1$.
As in \tough, let $(W_i)$ be a random walk whose increments
are standard exponential random variables, and write $U_{i,N}$ 
as $W_i/W_{N+1}$. Lemma 3.4.3 in \tough\ shows that 
$$
  \max_{\Lambda\leq n\leq N} \sharp\ROne=O_P(m_N\log N)
$$
as $\Lambda$ tends to infinity.

Robbins (1954) implies that provided $c$ is small enough, the set
$$
  \Omega=\{\, U_{i,N} \geq c i/N\, :\, 1\leq i\leq N\,\}
$$
has probability at least $1-\eta$. An integer $i$ in $\ROne$ is such that 
$U_{i,N}\leq m_n/n$, and on $\Omega$ we obtain $i\leq c m_n N/n$. So,
if $i$ is in $\ROne\setminus\{\, i^*\,\}$ and $\Omega$ occurs,
$$
  |V_{i^*}-V_i| \geq \min_{2\leq j\leq c m_n N/n} 
  V_{j,\lfloor cm_nN/n\rfloor}-V_{j-1,\lfloor cm_nN/n\rfloor}
  \, .
  \eqno{\equa{approxOneA}}
$$
Theorem 3.1 in Devroye (1982) implies that the right hand side of \approxOneA\
is almost surely at least $n^2/c m_n^2N^2\log N$ whenever $\Lambda$ is large
enough. For $n$ at least $\Lambda$, using that $n/m_n$ is asymptotically 
equivalent to a nondecreasing function and hence at least $\Lambda/m_\Lambda$, 
the right hand side of \approxOneA\ is
at least $c\Lambda^{-2(\beta+\epsilon)}/\log \Lambda$, and, if $\beta$ 
and $\epsilon$ are small enough, dominates $1/\sqrt N$ asymptotically. Since
$G_N=\Id+O_P(1/\sqrt{N})$ by Donsker's (1952) invariance principle,
$$
  \min_{i\in \ROne\setminus\{i^*\}} N|G_N(V_i)-G_N(V_{i^*})|
  \gsim {c n^2\over m_n^2 N\log N}
$$
for all $n$ in $(\Lambda,N)$, with probability at least $1-\eta$. Thus, writing
$n-NG_N(V_i)$ as $n-NG_N(V_{i^*})+N\bigl( G_N(V_{i^*})-G_N(V_i)\bigr)$, it
follows that
$$
  \max_{i\in\ROne\setminus\{\, i^*\,\}} g_{n-NG_N(V_i)}
  \lsim c g_{n^2/(m_n^2N\log N)}
  \leq c g_{\Lambda^2/\ml^2 N\log N} \, .
$$
Since $\prob{U_{1,N}\geq c/N}\geq 1-\eta$ if $c$ is small enough and $N$ is
large enough,
$$
  F^\leftarrow(1-U_{1,N})\leq c F^\leftarrow(1-1/N)
$$
with probability at least $1-\eta$. Using Potter's bound, this is at most
$c F^\leftarrow(1-1/n)(N/n)^{(1/\alpha)+\eta}$. 

Thus, with probability at least $1-\eta$, \approxOneB\ is at most
$$
  c m_N (\log N) g_{\Lambda^2\over m_n^2 N\log N} 
  \Bigl({N\over n}\Bigr)^{(1/\alpha)+\eta}F^\leftarrow(1-1/n) \, .
$$
For this bound to exceed $T\g n/\kl$ we must have, as $\Lambda$ tends to
infinity,
$$
  c m_N (\log N) g_{\Lambda^2\over \ml^2 N\log N}
  \Bigl({N\over\Lambda}\Bigr)^{(1/\alpha)+\eta} 
  \geq T {k(n)\over \kl}
  \gsim T \, .
  \eqno{\equa{approxOneC}}
$$
The left hand side is a regularly varying of $\Lambda$, of index
$$
  \beta(1+\epsilon) +(1-2\beta-\epsilon)(\gamma-1) 
  +\epsilon\Bigl({1\over\alpha}
  +\eta\Bigr)
  =\gamma-1+O(\beta)+O(\epsilon) \, .
$$
This index is negative if $\beta$ and $\epsilon$ are small enough, 
and \approxOneC\ cannot hold. This proves the lemma.\hfill\qed

\bigskip

Note that by construction $T_{1,n,N}^+$ is an approximation of the 
sum $\sum_{0\leq i<n} g_i X_{n-i}\One\{\, X_{n-i}>b_n\,\}$ by a single one
of its summands. Since each summand is at most $|g|_\infty X_{n,n}$, we see
that in order to show that
$$
  \limT \limsupL \prob{\existsnl T_{1,n,N}^+>T\g n/\kl} = 0 \, ,
$$
it suffices to prove that
$$
  \limT\limsupL \prob{ \existsnl X_{n,n}>T \g n/\kl} = 0 \, .
$$
Writing $X_{n,n}=F^\leftarrow(1-U_{1,n})$, this amounts to proving that
$$\displaylines{\quad
  \limT\limsupL \Prob{\existsnl
  \hfill\cr\hfill
  U_{1,n}\leq \oF\Bigl( T{k(n)\over\kl}F^\leftarrow(1-1/n)\Bigr)} = 0 \, .
  \quad\equa{eqB}\cr}
$$
Let $(V_i)$ be a new sequence of independent random variables having a uniform
distribution over $[\,0,1\,]$. Write $(V_{i,n})_{1\leq i\leq n}$ for the
order statistics of $(V_i)_{i\leq i\leq n}$. Setting $V_{1,0}=1$, we 
have $(U_{1,n})_{n\geq\Lambda}\eqd 
(U_{1,\Lambda}\wedge V_{1,n-\Lambda})_{n\geq\Lambda}$. Applying Bonferroni's
inequality, we see that for \eqB\ to hold it suffices to have
$$\displaylines{\quad
  \limT \limsupL \Prob{\existsnl 
  \hfill\cr\hfill
  U_{1,\Lambda}\leq \oF\Bigl(T{k(n)\over\kl}
  F^\leftarrow(1-1/n)\Bigr)} = 0
  \quad\equa{eqC}\cr}
$$
and, replacing $n$ by $\Lambda+n$, and setting
$$
  v_n
  =\oF\Bigl(T{k(\Lambda+n)\over\kl}F^\leftarrow\Bigl(1-{1\over \Lambda+n}\Bigr)
      \Bigr)
$$
to also have
$$
  \limT\limsupL \Prob{ \exists n\in (0,\Lambda^{1+\epsilon}) \, : \, 
  V_{1,n}\leq v_n} = 0 \, .
  \eqno{\equa{eqD}}
$$
The right hand side of the inequality involved in \eqC\ is equivalent to a 
function decreasing in $n$. So in the range of $n$ between $\Lambda$ 
and $\Lambda^{1+\epsilon}$ it is at most a constant times
$$
  T^{-\alpha} \oF\Bigl( F^\leftarrow\Bigl(1-{1\over\Lambda}\Bigr)\Bigr)
  \sim {T^{-\alpha}\over\Lambda} \, .
$$
Since the distribution of $\Lambda U_{1,\Lambda}$ converges to a standard
exponential one, \eqC\ holds.

To prove that \eqD\ holds, we use a blocking argument. Consider a real 
number $\theta$ greater than $1/(\alpha\gamma-1)$ and for any integer $p$ set 
$n_p=\lfloor \Lambda p^\theta\rfloor$. Potter's bound implies as $\Lambda$
tends to infinity and uniformly in $p$ positive
$$\eqalign{
  v_{n_p}
  &{}\lsim T^{-\alpha} \oF\Bigl( (1+p^\theta)^{\gamma-(1/\alpha)-\eta}
             F^\leftarrow\Bigl(1-{1\over\Lambda (1+p^\theta)}\Bigr)\Bigr) \cr
  &{}\lsim T^{-\alpha} \oF\Bigl( (1+p^\theta)^{\gamma-2\eta} 
             F^\leftarrow\Bigl(1-{1\over\Lambda}\Bigr)\Bigr)\cr
  &{}\lsim T^{-\alpha} (1+p^\theta)^{-\alpha(\gamma-3\eta)} {1\over\Lambda} 
             \, .\cr}
$$
Therefore, for any $\Lambda$ large enough and any positive $p$,
$$\eqalignno{
  \prob{ V_{1,n_p}\leq v_{n_p} }
  &{}\leq{}\Prob{ V_{1,n_p}\leq c T^{-\alpha} 
           (1+p)^{-\theta\alpha(\gamma-3\eta)} {1\over\Lambda}} \cr
  &{}\leq 1-\Bigl( 1-c T^{-\alpha}(1+p)^{-\theta\alpha(\gamma-3\eta)}
            {1\over\Lambda}\Bigr)^{\lfloor \Lambda p^\theta\rfloor} .
   \qquad\quad
  &\equa{eqE}\cr}
$$
Note that $(1+p)^{-\theta\alpha (\gamma-3\eta)}/\Lambda$ tends to $0$ as
$\Lambda$ tends to infinity, uniformly in $p$ nonnegative. So \eqE\ is at most
$$
  c \Lambda p^\theta T^{-\alpha}(1+p)^{-\theta\alpha(\gamma-3\eta)}
  {1\over\Lambda}
  \lsim c T^{-\alpha} p^{\theta(1-\alpha\gamma+3\alpha\eta)} \, .
$$
Given our choice of $\theta$, we see that if $\eta$ is small enough, the
exponent of $p$ is less than $-1$. Thus, Bonferroni's inequality implies
$$
  \limT\limsupL \prob{ \exists p\geq 1 \, : \,  V_{1,n_p}\leq v_{n_p} }
  = 0 \, .
$$
If $n$ is between $n_{p-1}$ and $n_p$, then $V_{1,n}\geq V_{1,n_p}$ and since
$v_{n_p}/v_{n_{p-1}}$ is bounded away from $0$ and infinity unformly in $p$
as $\Lambda$ tends to infinity, we proved \eqD\ and \eqB\ as well.\hfill\qed

\bigskip

\noindent{\it Step 4.} Let $T_n^-=\sum_{0\leq i<n} g_i X_{n-i}
\One_{(-\infty,a_n)}(X_{n-i})$. Combining steps 1, 2 and 3, we see that
$$\displaylines{
  \limT\limsupL \prob{\exists n>\Lambda \,:\, S_n-T_n^--\E (S_n-T_n^-)
  > T \g n/\kl} 
  \hfill\cr\hfill
  = 0 \, .\cr}
$$
Hence, in order to prove \eqKey, it suffices to show that
$$
  \limT\limsupL \prob{\exists n>\Lambda\,:\, |T_n^--\E T_n^-|>T\g n/\kl}
  = 0 \, .
$$
This follows by the very same arguments as in section 3.5 of \tough.

\bigskip

It remains to show that the limiting random variable involved in Theorem 
\mainTh.ii is almost surely finite. We write $\nu_+$ and $\nu_-$ for the 
restriction of $\nu$ to $(-\infty,0)$ and $(0,\infty)$ respectively. Let $N_+$
and $N_-$ be two independent Poisson processes with respective mean measures
$\lambda\otimes \nu_-$ and $\lambda\otimes \nu_+$. For a point 
process $N=\sum_{i\geq 1}\delta_{(t_i,x_i)}$ write $N^g$ 
for $\sup_{i\geq 1\atop j\geq 0} g_j x_i-t_i^\gamma$. Since $N_-$ 
and $N_+$ are independent, $N_-+N_+$ is a Poisson process with mean
intensity $\lambda\otimes\nu$. Since $(N_-+N_+)^g=N_+^g\vee N_-^g$, it suffices
to show that $N_+^g$ is finite. Write $N_+=\sum_{i\geq 1}\delta_{(t_i,x_i)}$. 
Since $(g_j)$ is bounded,
$$
  N_+^g\leq \sup_{i\geq 1} c x_i -t_i^\gamma
  \eqno{\equa{nPlusG}}
$$
whenever $c$ is at least $\max_{j\geq 0}g_j$. So it suffices to show that the
upper bound in \nPlusG\ is almost surely finite.

Since $N_+$ is a Poisson process, the random variables
$$
  M_k=\sup_{i:t_i\in [k,k+1)} cx_i -k^\gamma \, , \qquad k\geq 0\, ,
$$
are independent. Moreover,
$$
  \sup_{i\geq 1} c x_i-t_i^\gamma \leq \sup_{k\geq 0} M_k \, .
$$
Recall that $N_+$ has intensity $\lambda\otimes\nu_+$. Since
$$\eqalign{
  \prob{M_k>x}
  &{}=\Prob{ \exists i\,:\, (t_i,x_i)\in [k,k+1)\times
    \Bigl({x+k^\gamma\over c},\infty\Bigr)} \cr
  &{}=\Prob{ N\Bigl( [k,k+1)\times
    \Bigl({x+k^\gamma\over c},\infty\Bigr) \Bigr) \geq 1} \cr
  &{}=1-\exp\Bigl( -p\Bigl( {x+k^\gamma\over c}\Bigr)^{-\alpha}\Bigr) \, , \cr
  }
$$
we have
$$
  \prob{\max_{k\geq 1}M_k>x}
  \leq\sum_{k\geq 0}\biggl( 1-\exp\Bigl(-p\Bigl({c\over x+k^\gamma}\Bigr)^\alpha
          \Bigr)\biggr)
  \eqno{\equa{conclusionB}}
$$
This series is convergent since its $k$-th term is equivalent to
$c/k^{\alpha\gamma}$
as $k$ tends to infinity and $\alpha\gamma$ is greater than $1$. Bounding
$c/(x+k^\gamma)$ by $c/(1+k^\gamma)$ when $x$ exceeds $1$, the dominated 
convergence theorem implies that \conclusionB\ tends to $0$ as $x$ tends
to infinity, concluding the proof of Theorem \mainTh.

\bigskip


\noindent {\bf References.}

\medskip

{\leftskip=\parindent \parindent=-\parindent
 \par

Ph.\ Barbe, W.P.\ McCormick (2010a). Invariance principles for some FARIMA
and nonstationary linear processes in the domain of attraction of a stable
distribution, {\tt arXiv:1007.0576}.

Ph.\ Barbe, W.P.\ McCormick (2010b). Ruin probabilities in tough times --
Part 1: heavy-traffic approximation for fractionally integrated random 
walks in the domain of attraction of a nonGaussian stable distribution, 
{\tt arXiv:1101.4437}.

N.H.\ Bingham, C.M.\ Goldie, J.L.\ Teugels (1989). {\sl Regular Variation},
2nd ed. Cambridge University Press.

M.\ Donsker (1952). Justification and extension of Doob's heuristic approach
to the Kolmogorov-Smirnov theorems, {\sl Ann.\ Math.\ Statist.}, 23, 277--281.

J.\ Kiefer (1972). Iterated logarithm analogues for sample quantiles when 
$p_n\downarrow 0$, {\sl Proc.\ Sixth Berkeley Sympos.\ on Math.\ Statist.\ and
Probab.}, 1, 227--244.

H.\ Robbins (1954). A one-sided confidence interval for an unknown distribution
function, {\sl Ann.\ Math.\ Statist.}, 25, 409.

G.R.\ Shorack, J.A.\ Wellner (1978). Linear bounds on the empirical 
distribution function, {\sl Ann.\ Probab.}, 6, 349--353.

}

\bigskip

\setbox1=\vbox{\halign{#\hfil&\hskip 40pt #\hfill\cr
  Ph.\ Barbe            & W.P.\ McCormick\cr
  90 rue de Vaugirard   & Dept.\ of Statistics \cr
  75006 PARIS           & University of Georgia \cr
  FRANCE                & Athens, GA 30602 \cr
  philippe.barbe@math.cnrs.fr                        & USA \cr
                        & bill@stat.uga.edu \cr}}%
\box1%

\bye

mpage -2 -c -o -M-100rl-80b-220t -t veraverbeke7.ps > toto.ps

\bye